\documentclass[10pt,a4paper]{article}
\usepackage{amsmath,amssymb}
\usepackage[all]{xy}
\author{Fr{\'e}d{\'e}ric Chapoton}

\title{Op{\'e}rades diff{\'e}rentielles gradu{\'e}es sur les simplexes et les
  permuto{\`e}dres}

\newcommand{\ZZ}{\mathbb{Z}}
\newcommand{\QQ}{\mathbb{Q}}

\newcommand{\sym}{\mathfrak{S}}
\newcommand{\ext}{\operatorname{Ex}}
\newcommand{\card}{\operatorname{Card}}
\newcommand{\courb}{\curvearrowleft}
\newcommand{\id}{\operatorname{Id}}
\newcommand{\NS}{\not\Sigma}

\newcommand{\zin}{\operatorname{Zin}}
\newcommand{\dend}{\operatorname{Dend}}
\newcommand{\prelie}{\operatorname{PreLie}}

\newcommand{\com}{\operatorname{Com}}
\newcommand{\as}{\operatorname{As}}
\newcommand{\lie}{\operatorname{Lie}}

\newcommand{\pasc}{\operatorname{Pasc}}
\newcommand{\trias}{\operatorname{Trias}}
\newcommand{\pidual}{\amalg}

\newcommand{\perm}{\operatorname{Perm}}
\newcommand{\dias}{\operatorname{Dias}}
\newcommand{\leib}{\operatorname{Leib}}

\newcommand{\PP}{\mathcal{P}}

\newcommand{\FF}{\mathcal{F}}
\newcommand{\quot}{\operatorname{Quot}}
\newcommand{\rel}{\operatorname{Rel}}

\newtheorem{rema}{Remarque}

\newtheorem{defi}{D{\'e}finition}
\newtheorem{prop}{Proposition}

\newenvironment{proof}{{\sc Preuve : ~}}{\hfill\ensuremath{\square}}
\date{}
\begin{document}





\maketitle

\begin{abstract}
	We define several differential graded operads, some of them
	being related to families of polytopes : simplices and
	permutohedra. We also obtain a presentation by generators and
	relations of the operad $K$ on associahedra introduced in a
	previous article.
\end{abstract}

\section*{Introduction}

Les alg{\`e}bres associatives commutatives, les alg{\`e}bres associatives et
les alg{\`e}bres de Lie font partie du paysage math{\'e}matique depuis
longtemps.

Plus r{\'e}cemment sont apparues les alg{\`e}bres de Leibniz, les dig{\`e}bres,
les alg{\`e}bres dendriformes et les alg{\`e}bres zinbiel, toutes introduites
par Loday avec des motivations en K-th{\'e}orie, voir
\cite{loday-enseign,loday-cras,loday-overview,dialgebra}.

Ces types d'alg{\`e}bres s'organisent en un diagramme d'op{\'e}rades,
qui se compl\-{\`e}te par les alg{\`e}bres permutatives et les
alg{\`e}bres pr{\'e}-Lie
\footnote{aussi connues sous le nom d'alg{\`e}bres sym{\'e}triques {\`a} gauche}, voir
\cite{endo,rooted}.

Par ailleurs, la relation entre l'op{\'e}rade des alg{\`e}bres dendriformes et
les sommets des polytopes de Stasheff m{\`e}ne {\`a} la d{\'e}finition d'une
op{\'e}rade diff{\'e}rentielle gradu{\'e}e $K$ sur toutes les faces des polytopes
de Stasheff \cite{trans}. La place que prend cette op{\'e}rade parmi les
pr{\'e}c{\'e}dentes incite {\`a} compl{\'e}ter encore le diagramme, ce que nous
faisons ici.

On va d{\'e}finir les seconde et quatri{\`e}me lignes du diagramme commutatif
suivant, o{\`u} les fl{\`e}ches marqu{\'e}es \textit{q.i.} sont des
quasi-isomorphismes.

\begin{equation}
  \label{grosdiag}
  \xymatrix{
    \relax\zin & \dend \ar@{->>}[l] & \prelie \ar@{->}[l] \\
    \Pi \ar@{->>}[u] & K \ar@{->>}[l] \ar@{->>}[u] & \Lambda \ar@{->}[l] \ar@{->>}[u] \\
    \com \ar@{->}[u]^{q.i.} & \as \ar@{->>}[l]
    \ar@{->}[u]^{q.i.} & \lie \ar@{->}[l] \ar@{->}[u]^{q.i.} \\
    \pasc \ar@{->>}[u]^{q.i.} & \trias \ar@{->>}[l]
    \ar@{->>}[u]^{q.i.} &
    \pidual \ar@{->}[l] \ar@{->>}[u]^{q.i.}\\
    \perm \ar@{->}[u] & \dias \ar@{->>}[l] \ar@{->}[u] & \leib \ar@{->}[l] \ar@{->}[u] }
\end{equation}

La sym{\'e}trie centrale de ce diagramme correspond {\`a}
la dualit{\'e} de Koszul des op{\'e}rades quadratiques binaires.

\bigskip

L'op{\'e}rade $K$ a {\'e}t{\'e} d{\'e}finie dans la premi{\`e}re partie de \cite{trans} ;
elle est en relation avec les associa{\`e}dres. Ici, on commence par
d{\'e}finir deux op{\'e}rades $\Pi$ (faisant intervenir les permuto{\`e}dres) et
$\pasc$ (li{\'e}e aux simplexes). On donne ensuite une pr{\'e}sentation par
g{\'e}n{\'e}rateurs et relations de ces op{\'e}rades, qui montre en particulier
que ces trois op{\'e}rades sont quadratiques binaires. Les r{\'e}sultats de
Markl \cite{markl96} sur les lois distributives permettent de
d{\'e}montrer que $\Pi$ est de Koszul.

On donne ensuite une pr{\'e}sentation des op{\'e}rades $\pidual$, $\Lambda$ et
$\trias$, duales quadratiques des pr{\'e}c{\'e}dentes au sens de Ginzburg et
Kapranov \cite{ginzburg-kapranov}. L'op{\'e}rade $\trias$ provient d'une
op{\'e}rade non-sym{\'e}\-trique $\trias'$. On d{\'e}montre que $\trias'$ est
isomorphe {\`a} l'op{\'e}rade $\pasc$ vue comme op{\'e}rade non-sym{\'e}trique.

Durant la pr{\'e}paration de cet article, j'ai re{\c c}u une note de Loday et
Ronco \cite{loday-annonce} annon{\c c}ant des r{\'e}sultats similaires pour des
variantes non-gradu{\'e}es des op{\'e}rades $K$ et $\trias$.


\section{Op{\'e}rade $\Pi$}

On se place dans toute la suite de l'article dans la cat{\'e}gorie
mono{\"\i}dale sym{\'e}trique des complexes de $\ZZ$-modules.

On d{\'e}finit une op{\'e}rade diff{\'e}rentielle gradu{\'e}e $\Pi$. Les complexes
sous-jacents sont les complexes de cocha{\^\i}nes cellulaires des
permuto{\`e}dres.

\subsection{Complexes sous-jacents}

Soit $I$ un ensemble fini. On note $\Pi(I)$ le $\ZZ$-module libre
engendr{\'e} par les expressions de la forme $ \pi_1\otimes \pi_2 \otimes
\dots \otimes \pi_p$, o{\`u} $\pi_j=i_{j,1}\wedge \dots \wedge i_{j,p_j}$
avec $i_{j,k}\in I$ et o{\`u} chaque {\'e}l{\'e}ment de $I$ appara{\^\i}t dans un
certain $\pi_j$ et un seul, modulo les relations d'antisym{\'e}trie des
produits ext{\'e}rieurs. Ces expressions correspondent {\`a} des partitions
ordonn{\'e}es de $I$, en bijection avec les faces du permuto{\`e}dre
correspondant {\`a} $I$.

Si $\pi_j=i_{j,1}\wedge \dots \wedge i_{j,p_j}$, on appelle dimension
de $\pi_j$, not{\'e} $\dim(\pi_j)$, l'entier $p_j-1$. On munit $\Pi(I)$ de
la graduation par la dimension en posant, si $\pi= \pi_1\otimes \pi_2
\otimes \dots \otimes \pi_p$,
\begin{equation*}
  \dim(\pi)=\sum_{j=1}^p \dim(\pi_j),
\end{equation*}
et d'une diff{\'e}rentielle $d$ de dimension $+1$ d{\'e}finie par la formule
suivante :
\begin{equation}
  d(\pi)=\sum_{j=1}^{p-1}(-1)^{\dim(\pi_1)+\dots+\dim(\pi_j)}\pi_1 \otimes
  \dots\otimes \pi_j\wedge \pi_{j+1}\otimes \dots \otimes \pi_p.
\end{equation}

On v{\'e}rifie sans difficult{\'e} que $d$ est de carr{\'e} nul.

\subsection{Composition de l'op{\'e}rade}

Soient $I$ et $J$ deux ensembles finis et $i\in I$. Soient
$\pi=\pi_1\otimes \pi_2 \otimes \dots \otimes \pi_p\in\Pi(I)$ et
$\mu=\mu_1\otimes \mu_2 \otimes \dots \otimes \mu_q\in\Pi(J)$. On va
d{\'e}finir la composition $\pi\circ_i \mu$, appartenant {\`a}
$\Pi(I\setminus\{i\}\sqcup J)$. Soit $\ell$ l'unique indice tel que
$i$ appara{\^\i}t dans $\pi_{\ell}$. Quitte {\`a} changer de signe, on peut
supposer que $\pi_{\ell}=\pi'_{\ell} \wedge i$, avec la convention
que, si $\pi_{\ell}=i$, alors $\pi'_{\ell}\wedge \nu=\nu$ pour tout
$\nu$.

On pose alors
\begin{equation}
  \pi\circ_i \mu=\sum_{\sigma}\varepsilon(\sigma) \pi_1\otimes \dots \otimes 
  \pi_{\ell-1}\otimes\pi'_{\ell}\wedge
  \mu_1\otimes \sigma(\pi_{\ell+1}\otimes \dots \otimes \pi_p\,,\,\mu_2\otimes
  \dots \otimes \mu_q),
\end{equation}
o{\`u} $\sigma$ d{\'e}crit les battages de $\pi_{\ell+1}\,,\, \dots\,,\,
\pi_p$ avec $\mu_2\,,\,\dots \,,\, \mu_q$, $\sigma(\cdots)$ est le
produit tensoriel de $\pi_{\ell+1}\,,\, \dots\,,\,
\pi_p\,,\,\mu_2\,,\,\dots \,,\, \mu_q$ dans l'ordre sp{\'e}cifi{\'e} par
$\sigma$ et $\epsilon(\sigma)$ d{\'e}signe le signe obtenu par la r{\`e}gle de
Koszul (en utilisant la volte gradu{\'e}e par la dimension) pour passer de
$\pi_{\ell+1}\otimes \dots\otimes\pi_p\otimes\mu_1\otimes\dots\otimes
\mu_q$ {\`a} $\mu_1\otimes \sigma(\pi_{\ell+1}\otimes \dots \otimes
\pi_p\,,\,\mu_2\otimes \dots \otimes \mu_q)$.

\begin{prop}
  $\Pi$ est une op{\'e}rade diff{\'e}rentielle gradu{\'e}e.
\end{prop}
\begin{proof}
  L'{\'e}l{\'e}ment $1$ de $\Pi(\{1\})$ est une unit{\'e}. La preuve que la
  composition est un morphisme de complexes est un calcul long mais
  sans difficult{\'e}, de m{\^e}me que celle de l'associativit{\'e}.
  L'{\'e}quivariance, i.e la fonctorialit{\'e} de la d{\'e}finition sur le
  groupo{\"\i}de des ensembles finis, est imm{\'e}diate.
\end{proof}

\subsection{Pr{\'e}sentation par g{\'e}n{\'e}rateurs et relations}

\begin{prop}
  \label{pi-binaire}
  L'op{\'e}rade $\Pi$ est engendr{\'e}e par $\Pi(\{1,2\})$. 
\end{prop}
\begin{proof}
  On montre d'abord par r{\'e}currence que $1\wedge 2 \wedge \dots \wedge
  n$ s'obtient par composition it{\'e}r{\'e}e de $1\wedge2$ en utilisant la
  relation
  \begin{equation}
    (1\wedge 2 \wedge \dots \wedge n-1 \wedge
    \star)\circ_{\star}(n\wedge n+1)=(1\wedge 2 \wedge \dots \wedge n+1).
  \end{equation}
  On montre de m{\^e}me que $1\otimes 2\otimes \dots \otimes n$ s'obtient
  par composition it{\'e}r{\'e}e de $1 \otimes 2$ en utilisant la relation
  \begin{equation}
    (1\otimes 2 \otimes \dots \otimes n-1 \otimes
    \star)\circ_{\star}(n\otimes n+1)=(1\otimes 2 \otimes \dots \otimes n+1).
  \end{equation}
  En effectuant successivement la composition d'un {\'e}l{\'e}ment de la forme
  $i_1 \wedge i_2 \wedge \dots \wedge i_p$ dans chaque terme d'un
  {\'e}l{\'e}ment de la forme $1\otimes 2\otimes \dots \otimes q$, on obtient
  le r{\'e}sultat voulu.
\end{proof}

\smallskip

On a dans $\Pi$ les relations suivantes :
\begin{align}
  \label{rel1}
  (1\otimes 2)\circ_1 (1\otimes 2)&=(1\otimes 2)\circ_2 (1\otimes 2)+(1\otimes 2)\circ_2 (2\otimes 1),\\
  \label{rel2}
  (1\wedge 2)\circ_2(1\otimes 2)&=(1\otimes 2)\circ_1(1\wedge2)=\tau (1\wedge 2)\circ_1(2\otimes 1),\\
  \label{rel3}
  (2\wedge 1)\circ_1(1\wedge2)&=(1\wedge 2)\circ_2 (1\wedge2),
\end{align}
o{\`u} $\tau$ est le cycle $3\to 2\to 1\to 3$ agissant dans
$\Pi(\{1,2,3\})$.

\begin{prop}
  L'op{\'e}rade $\Pi$ est isomorphe au quotient de l'op{\'e}rade libre sur
  $\Pi(\{1,2\})$ par l'id{\'e}al engendr{\'e} par les relations ci-dessus. En
  particulier $\Pi$ est une op{\'e}rade quadratique binaire.
\end{prop}
\begin{proof}
  Soit $\FF$ l'op{\'e}rade libre sur $\Pi(\{1,2\})$, $\rel$ les relations
  ci-dessus et $\quot$ l'op{\'e}rade quotient de $\FF$ par l'id{\'e}al
  engendr{\'e} par $\rel$. On a un morphisme d'op{\'e}rades de $\quot \to
  \Pi$, surjectif par la proposition \ref{pi-binaire}.
  
  On remarque que la relation (\ref{rel1}) est celle de l'op{\'e}rade
  $\zin$ et que la relation (\ref{rel3}) est celle de la suspension de
  l'op{\'e}rade $\com$. Par ailleurs, on v{\'e}rifie que la relation
  (\ref{rel2}) d{\'e}finit une loi distributive au sens de Markl
  \cite{markl96} reliant l'op{\'e}rade $\zin$ et la suspension de $\com$.
  Comme les $\ZZ$-modules sous-jacents {\`a} $\zin$ et $\com$ sont libres,
  les $\ZZ$-modules sous-jacents {\`a} $\quot$ sont {\'e}galement libres.
  
  Pour montrer que la projection $\quot\to\Pi$ est un isomorphisme, il
  suffit de montrer l'{\'e}galit{\'e} des rangs des $\ZZ$-modules libres
  sous-jacents. On rappelle la d{\'e}finition de la s{\'e}rie g{\'e}n{\'e}ratrice
  d'une op{\'e}rade diff{\'e}rentielle gradu{\'e}e $\PP$ :
  \begin{equation*}
    g_{\PP}(x,t)=\sum_{n\geq 1}\sum_{k}\dim \PP^k(n) (-t)^k \frac{x^n}{n!}.
  \end{equation*}
  Par une propri{\'e}t{\'e} des lois distributives, la s{\'e}rie g{\'e}n{\'e}ratrice de
  $\quot$ s'obtient par composition de celles de $\zin$ et de la
  suspension de $\com$. On v{\'e}rifie qu'elle co{\"\i}ncide avec celle de
  $\Pi$, voir le tableau r{\'e}capitulatif {\`a} la fin de l'article.
\end{proof}
\smallskip

\medskip

La diff{\'e}rentielle de $\Pi$ est donc uniquement d{\'e}termin{\'e}e par le respect
de la structure d'op{\'e}rade et les conditions suivantes :
\begin{align}
  d(1\otimes 2)&=1\wedge 2,\\
  d(1\wedge 2)&=0.
\end{align}

On en d{\'e}duit la description suivante des $\Pi$-alg{\`e}bres.
\smallskip

\begin{defi}
  Une $\Pi$-alg{\`e}bre est la donn{\'e}e d'un complexe de cocha{\^\i}nes $V$ et
  d'applications $\prec \,: V\otimes V\to V$ et $\times : V\otimes
  V\to V[-1]$ telles que pour tous $x,y,z$ dans $V$,
  \begin{align} 
    (x\prec y)\prec z&=x\prec(y\prec z)+(-1)^{yz}x\prec(z\prec y),\\
    x\times (y\prec z)&=(x\times y)\prec z=(-1)^{yz+z}(x\prec z)\times y,\\
    x\times y&=(-1)^{x y+x+y+1}y\times x,\\
    (x\times y)\times z&=x\times (y\times z)
  \end{align}
  et 
  \begin{align}
    (dx\prec y)+(-1)^{x}(x\prec d y)-d (x\prec y)=(-1)^{x}x\times y,\\  
    (dx\times y)+(-1)^{x+1}(x\times d y)-d (x\times y)=0.  
  \end{align}
\end{defi}

\begin{rema} 
On utilise ici et plus loin l'abus de notation qui consiste, dans les
exposants de $(-1)$, {\`a} {\'e}crire $x,y,z$ pour $\dim(x),\dim(y),\dim(z)$.
\end{rema}

\begin{prop}
  L'op{\'e}rade $\Pi$ est de Koszul.
\end{prop}
\begin{proof}
  Soit $\Pi'$ l'op{\'e}rade obtenue en munissant $\Pi$ de la
  diff{\'e}rentielle nulle. Les op{\'e}rades $\com$ et $\zin$ sont de Koszul
  et $\Pi'$ est d{\'e}crite par une loi distributive, donc $\Pi'$ est de
  Koszul par le th{\'e}or{\`e}me 4.5 de \cite{markl96}. Par ailleurs, la
  koszulit{\'e} d'une op{\'e}rade diff{\'e}rentielle dont les complexes
  sous-jacents sont born{\'e}s est une cons{\'e}quence, par un argument de
  suite spectrale, de la koszulit{\'e} de la m{\^e}me op{\'e}rade avec la
  diff{\'e}rentielle nulle. Ceci entra{\^\i}ne que $\Pi$ est de Koszul.
\end{proof}

\smallskip

On observe, par ailleurs, que le quotient de $\Pi$ par l'id{\'e}al
diff{\'e}rentiel form{\'e} des {\'e}l{\'e}ments de dimension non nulle est une op{\'e}rade
(de diff{\'e}rentielle nulle) isomorphe {\`a} l'op{\'e}rade $\zin$ des alg{\`e}bres
zinbiel.

\smallskip

D'autre part, on a un morphisme de $\com$ muni de la diff{\'e}rentielle
nulle dans $\Pi$. C'est un quasi-isomorphisme car les complexes
sous-jacents sont les complexes cellulaires des permuto{\`e}dres, de sorte
que l'homologie est en degr{\'e} $0$ avec action triviale des groupes
sym{\'e}triques.

\section{Op{\'e}rade $\pasc$}

On d{\'e}finit une op{\'e}rade diff{\'e}rentielle gradu{\'e}e $\pasc$. Les complexes
sous-jacents sont les complexes de cha{\^\i}nes des simplexes.

\subsection{Complexes sous-jacents}

Si $I$ est un ensemble fini non vide, on note $\ZZ^I$ le $\ZZ$-module
libre de base $(e^I_i)_{i\in I}$. Soit $\ext(I)$ l'alg{\`e}bre ext{\'e}rieure
avec unit{\'e} sur $\ZZ^I$, munie de la graduation standard.

Soit $i\in I$, on note $ \theta_i^I$ la d{\'e}rivation ({\`a} gauche) de degr{\'e}
$-1$ de $\ext(I)$ d{\'e}finie sur les g{\'e}n{\'e}rateurs par
\begin{equation*}
  \theta^I_i(e^I_j)=\delta_{ij}.
\end{equation*}

On munit $\ext(I)$ de la diff{\'e}ren\-tielle $d$ de degr{\'e} $-1$ d{\'e}finie
par la formule 
\begin{equation*}
 d(x)=\sum_{i\in I}\theta^I_i(x).
\end{equation*}

On d{\'e}finit $\pasc(I)$ comme le complexe quotient de $\ext(I)$ par
l'unit{\'e} de $\ext(I)$, muni d'une graduation d{\'e}cal{\'e}e nomm{\'e}e dimension :
si $P$ est une partie de non vide de $I$, on pose
\begin{equation*}
  \dim\left(\pm \wedge_{i\in P}\,e^I_i\right)=\card(P)-1.
\end{equation*}

On note encore $\theta^I_i$ l'application obtenue par passage au
quotient. Elle v{\'e}rifie
\begin{equation}
  \theta^I_i(x\wedge y)=\theta^I_i(x)\wedge y+
  (-1)^{1+\dim(x)}x\wedge\theta^I_i(y),
\end{equation}
pour $x,y$ homog{\`e}nes de dimensions non nulles dans $\pasc(I)$. 


\subsection{Composition}

Soient $I,J$ deux ensembles finis, $P\subset I$ et $Q\subset J$ des
parties non vides et $i\in I$. On choisit une {\'e}num{\'e}ration arbitraire
de $P$ et $Q$ : $P=\{i_1,\dots,i_p\}$ et $Q=\{j_1,\dots,j_q\}$.

Soient $x=e^I_{i_1}\wedge\dots\wedge e^I_{i_p}$ et
$y=e^J_{j_1}\wedge\dots\wedge e^J_{j_q}$. On d{\'e}finit la composition $x
\circ_i y$, appartenant {\`a} $\pasc(I\setminus\{i\}\sqcup J)$, par
\begin{equation}
  x \circ_i y=
  \begin{cases}
    (-1)^{\dim (x)} \theta^{I\sqcup J}_i(x\wedge y)\text{ si }i\in P,\\
    x \text{ si }i\not\in P\text{ et }\dim y=0,\\
    0 \text{ si }i\not\in P\text{ et }\dim y>0.
  \end{cases}
\end{equation}

\begin{rema}
Dans cette d{\'e}finition, les m{\^e}mes produits ext{\'e}rieurs sont
implicitement consid{\'e}r{\'e}s dans des alg{\`e}bres ext{\'e}rieures diff{\'e}rentes,
via les identifications canoniques.
\end{rema}

\begin{prop}
  $\pasc$ est une op{\'e}rade diff{\'e}ren\-tielle gradu{\'e}e.
\end{prop}

\begin{proof}
  L'{\'e}l{\'e}ment $e^{\{1\}}_1$ de $\pasc(\{1\})$ est clairement une unit{\'e}.
  La preuve que la composition est un morphisme de complexes passe par
  une {\'e}tude des diff{\'e}rents cas possibles. De m{\^e}me, on montre
  l'associativit{\'e} et l'{\'e}quivariance de la composition en distinguant
  les diff{\'e}rentes situations.
\end{proof}

\subsection{Pr{\'e}sentation par g{\'e}n{\'e}rateurs et relations}

\begin{prop}
  \label{pasc-binaire}
  L'op{\'e}rade $\pasc$ est engendr{\'e}e par $\pasc(\{1,2\})$.   
\end{prop}
\begin{proof}
  On montre d'abord par r{\'e}currence que, pour tout ensemble $I$ et tout
  $i\in I$, on peut obtenir $e^I_i$. On utilise la relation
  \begin{equation*}
  e^{I\setminus\{j,k\}\sqcup\{\star\}}_{i}\circ_{\star}
  e^{\{j,k\}}_{j}=e^I_i, \end{equation*} pour $i,j,k\in I$ deux {\`a} deux
  distincts.

  On montre de m{\^e}me que, pour tout ensemble $I$, on
  peut obtenir $\wedge_{i\in I} e^I_i$. On utilise la formule
  \begin{equation*}
    \bigg{(} \wedge_{i\in I\setminus\{j,k\}}\,
    e^{I\setminus\{j,k\}\sqcup\{\star\}}_{i}\wedge 
    e^{I\setminus\{j,k\}\sqcup\{\star\}}_{\{\star\}}\bigg{)}\circ_{\star}
    \bigg{(}e^{\{j,k\}}_{j}\wedge e^{\{j,k\}}_{k}\bigg{)}=
    \wedge_{i\in I}\, e^I_i,
  \end{equation*}
  pour $i,j,k\in I$ deux {\`a} deux distincts.

  Enfin, on d{\'e}duit le r{\'e}sultat voulu des deux points pr{\'e}c{\'e}dents en
  utilisant la relation
  \begin{equation*}
    e^{I\setminus P \sqcup\{\star\}}_{\star}\circ_{\star}
   \bigg{(} \wedge_{i\in P}\,  e^{P}_{i}\bigg{)}=\wedge_{i\in P}\, e^I_i,
  \end{equation*}
  pour $P$ une partie non vide de $I$.
\end{proof}

\smallskip

On a dans $\pasc$ les relations suivantes :
\begin{align}
  \label{rela1}
  e_1 \circ_1 e_1&=\tau^2 \, e_2\circ_1 e_1,\\
  \label{rela2}
  e_2 \circ_1 e_1&= e_2 \circ_1 e_2,\\
  \label{rela3}
  e_1\circ_1 e_2&=\tau\, e_2\circ_1 e_2,\\
  \label{rela4}
  e_2\, \circ_1 (e_1\wedge e_2)&=0,\\
  \label{rela5}
  \tau (e_1\wedge e_2) \circ_1 e_2 &= e_1 \circ_1 (e_1\wedge e_2),\\  
  \label{rela6}
  -\tau^2 (e_1\wedge e_2)\circ_1 e_1  &= e_1 \circ_1 (e_1\wedge e_2),\\
  \label{rela7}
  (e_1\wedge e_2)\circ_1 (e_1\wedge e_2)
  &=\tau(e_1\wedge e_2)\circ_1(e_1\wedge e_2),
\end{align}
o{\`u} $\tau$ est le cycle $3\to 2\to 1\to 3$.

\begin{prop}
  \label{pasc-pres}
  L'op{\'e}rade $\pasc$ est isomorphe au quotient de l'op{\'e}rade libre sur
  $\pasc(\{1,2\})$ par l'id{\'e}al engendr{\'e} par les relations ci-dessus.
  En particulier $\pasc$ est une op{\'e}rade quadratique binaire.
\end{prop}
\begin{proof}
  Soit $\FF$ l'op{\'e}rade libre sur $\pasc(\{1,2\})$, $\rel$ les
  relations ci-dessus et $\quot$ l'op{\'e}rade quotient de $\FF$ par
  l'id{\'e}al engendr{\'e} par $\rel$. Par la proposition \ref{pasc-binaire},
  on a un morphisme surjectif d'op{\'e}rades $\phi$ de $\quot \to \pasc$.
  La d{\'e}monstration consiste {\`a} construire un inverse $\psi$ par
  r{\'e}currence. On peut n{\'e}gliger les questions de signes et donc
  raisonner sur des parties plut{\^o}t que sur des produits ext{\'e}rieurs. On
  dit que la composition $P\circ_i Q$ est interne si $i\in P$ et
  externe sinon.
  
  Supposons donc un inverse $\psi$ construit jusqu'au rang $n-1$. Soit
  $I$ un ensemble {\`a} $n\geq 3$ {\'e}l{\'e}ments et $P$ une partie non vide de
  $I$. On distingue deux cas : ou bien $P$ contient au moins deux
  {\'e}l{\'e}ments, ou bien son compl{\'e}mentaire contient au moins deux
  {\'e}l{\'e}ments. Dans le premier cas, on peut {\'e}crire $P$ comme compos{\'e}
  interne de la partie $p=\{1,2\}$ de $\{1,2\}$ dans une partie $P'$
  d'un ensemble de cardinal $n-1$. Dans le second cas, on peut {\'e}crire
  $P$ comme compos{\'e} externe de la partie $p=\{1\}$ de $\{1,2\}$ dans
  une partie $P'$ d'un ensemble de cardinal $n-1$. Dans les deux cas,
  on pose alors $\psi(P)=\psi(P')\circ \psi(p)$. On a $\phi \psi=\id$
  car $\phi$ est un morphisme d'op{\'e}rades.
  
  Pour montrer que $\psi$ ne d{\'e}pend pas des choix faits, on fixe deux
  choix diff{\'e}rents et on distingue deux cas. Si les deux choix sont
  disjoints, au sens o{\`u} les paires d'{\'e}l{\'e}ments choisis sont sans
  intersection, on utilise l'axiome d'associativit{\'e} des op{\'e}rades. Si
  les deux choix ont un {\'e}l{\'e}ment en commun, on utilise les relations
  (\ref{rela1}), (\ref{rela2}), (\ref{rela3}) ou (\ref{rela7}).
  
  Les autres relations de $\pasc$ servent {\`a} choisir des repr{\'e}sentants
  particuliers des {\'e}l{\'e}ments de $\quot$ dans $\FF$. Plus pr{\'e}cis{\'e}ment,
  on {\'e}limine les compositions du type du cot{\'e} gauche des relations
  (\ref{rela4}), (\ref{rela5}) et (\ref{rela6}). Ces repr{\'e}sentants
  peuvent s'{\'e}crire comme composition de telle fa{\c c}on que leur image par
  $\phi$ soit une composition du type de celles utilis{\'e}es pour d{\'e}finir
  $\psi$.
 
  Ceci permet de d{\'e}montrer, en utilisant le fait que $\psi$ ne d{\'e}pend
  pas des choix, que $\psi \phi =\id$, donc $\psi$ est un inverse de
  $\phi$ pour les parties des ensembles {\`a} $n$ {\'e}l{\'e}ments, ce qui termine
  la r{\'e}currence.
\end{proof}

\medskip

La diff{\'e}rentielle de $\pasc$ est donc uniquement d{\'e}termin{\'e}e par le respect
de la structure d'op{\'e}rade et les conditions suivantes :
\begin{align}
  d(e_1\wedge e_2)&=e_1-e_2,\\
  d(e_1)&=0.
\end{align}

On en d{\'e}duit la description suivante des $\pasc$-alg{\`e}bres.
\smallskip

\begin{defi}
  Une alg{\`e}bre pascale est la donn{\'e}e d'un complexe de cha{\^\i}nes $V$,
  d'applications $ \dashv\, : V\otimes V\to V $ et $\times : V\otimes
  V \to V[-1]$ telles que, pour tous $x,y,z$ dans $V$,
  \begin{align}
    x\dashv (y\dashv z)&=(-1)^{y z}x\dashv (z \dashv y)=(x\dashv y)\dashv z,\\
    (x\times y)\dashv z&=x\times (y\dashv z),\\
    (x\dashv y)\times z&=(-1)^{yz+y}x\times (z\dashv y),\\
    x\dashv (y\times z)&=0,\\
    x\times y&=(-1)^{x y+x+y+1}y\times x,\\
    (x\times y)\times z&=x\times (y\times z)
  \end{align}
  et 
  \begin{align}
    (dx\times y)+(-1)^{x+1}(x\times dy)-d(x\times y)&=
    x\dashv y-(-1)^{xy} y\dashv x,\\
    (dx\dashv y)+(-1)^{x}(x\dashv dy)-d(x\dashv y)&=0.
  \end{align}
\end{defi}

\begin{rema}
On peut penser que $\pasc$ est de Koszul.
\end{rema}


\smallskip

On observe que les {\'e}l{\'e}ments de dimension nulle de $\pasc$ forment une
sous-op{\'e}rade (de diff{\'e}rentielle nulle) isomorphe {\`a} l'op{\'e}rade $\perm$
des dig{\`e}bres commutatives.

\smallskip

D'autre part, on a un morphisme quotient de $\pasc$ dans $\com$ muni
de la diff{\'e}rentielle nulle. C'est un quasi-isomorphisme car les
complexes sous-jacents de $\pasc$ sont les complexes simpliciaux de
simplexes, de sorte que l'homologie est en degr{\'e} $0$ avec action
triviale des groupes sym{\'e}triques.

\section{Op{\'e}rade $K$}

\subsection{Alg{\`e}bres sur associa{\`e}dres}

Dans un article ant{\'e}rieur \cite{trans}, on a d{\'e}fini une op{\'e}rade
diff{\'e}rentielle gradu{\'e}e $K$ sur les complexes de cocha{\^\i}nes cellulaires
des associa{\`e}dres. Cette op{\'e}rade provient d'une op{\'e}rade
non-sym{\'e}\-trique $K'$, qui est celle que l'on consid{\`e}re ici. On
travaille dans cette section avec des op{\'e}rades non-sym{\'e}\-triques. Tous
les r{\'e}sultats {\'e}nonc{\'e}s ci-dessous pour des op{\'e}rades non-sym{\'e}\-triques
peuvent {\^e}tre traduits en r{\'e}sultats similaires pour les op{\'e}rades
sym{\'e}triques associ{\'e}es.

\subsection{Pr{\'e}sentation par g{\'e}n{\'e}rateurs et relations}

On rappelle que $K'(n)$ est le complexe de cocha{\^\i}nes cellulaires de
l'associa{\`e}dre ou polytope de Stasheff de dimension $n-1$ dont les
cellules sont en bijection avec les arbres plans {\`a} $n+1$ feuilles. On
renvoie {\`a} \cite{trans} pour la d{\'e}finition de la composition.

\begin{prop}
  \label{K-binaire}
  L'op{\'e}rade $K'$ est engendr{\'e}e par $K'(\{1,2\})$. 
\end{prop}
\begin{proof}
  Par r{\'e}currence sur $n$. Soit donc $n\geq 3$ et $T$ un arbre dans
  $K'(n)$, la racine en bas. On consid{\`e}re un sommet maximal $s$ de $T$, au
  sens o{\`u} aucun sommet ne se trouve greff{\'e} au dessus de lui. 
  
  Si $s$ est binaire, on choisit de plus un cot{\'e}, droite ou gauche, et
  on obtient $T$ comme compos{\'e} d'un arbre binaire {\`a} trois feuilles
  dans un arbre $T'$ {\`a} $n$ feuilles obtenu en supprimant $s$.
  
  Si $s$ n'est pas binaire, on choisit de plus une feuille non-extr{\^e}me
  de ce sommet, et on obtient $T$ comme compos{\'e} d'une corolle {\`a} trois
  feuilles dans un arbre $T'$ {\`a} $n$ feuilles obtenu en supprimant
  l'ar{\^e}te issue de la feuille non-extr{\^e}me choisie.
  
  Comme $T'$ est compos{\'e} {\`a} partir d'{\'e}l{\'e}ments de $K'(\{1,2\})$ par
  hypoth{\`e}se de r{\'e}currence, il en va de m{\^e}me pour $T$. Ceci termine la
  r{\'e}currence.
\end{proof}
\smallskip

On a dans $K'$ les relations suivantes :
\begin{align}
  \label{relK1}
  (1 > 2)\circ_2 (1 < 2)&=(1 < 2)\circ_1 (1 > 2),\\
  \label{relK2}
  (1 > 2)\circ_2 (1 > 2)&=(1 > 2)\circ_1 (1 > 2)+(1 > 2)\circ_1 (1 < 2),\\
  \label{relK3}
  (1 < 2)\circ_1 (1 < 2)&=(1 < 2)\circ_2 (1 < 2)+(1 < 2)\circ_2 (1 > 2),\\
  \label{relK4}
  (1 > 2)\circ_2 (1 | 2)&= (1 | 2)\circ_1 (1 > 2),\\
  \label{relK5}
  (1 < 2)\circ_1 (1 | 2) &= (1 | 2)\circ_2 (1 < 2),\\
  \label{relK6}
  (1 | 2)\circ_1 (1 < 2)&=(1 | 2)\circ_2 (1 > 2),\\
  \label{relK7}
  (1 | 2)\circ_1 (1 | 2)&=-(1 | 2)\circ_2 (1 | 2).
\end{align}

\begin{prop}
  L'op{\'e}rade $K'$ est isomorphe au quotient de l'op{\'e}rade
  non-sym{\'e}\-trique libre sur $K'(\{1,2\})$ par l'id{\'e}al engendr{\'e} par
  les relations ci-dessus. En particulier $K'$ est une op{\'e}rade
  quadratique binaire.
\end{prop}
\begin{proof}
  Soit $\FF$ l'op{\'e}rade libre sur $K'(\{1,2\})$, $\rel$ les relations
  ci-dessus et $\quot$ l'op{\'e}rade quotient de $\FF$ par l'id{\'e}al
  engendr{\'e} par $\rel$. Par la proposition \ref{K-binaire}, on a un
  morphisme surjectif d'op{\'e}rades $\phi$ de $\quot \to K$. La
  d{\'e}monstration consiste {\`a} construire un inverse $\psi$ par
  r{\'e}currence.
  
  On a un inverse pour les arbres {\`a} $2$ feuilles. Supposons donc un
  inverse $\psi$ construit jusqu'aux arbres {\`a} $n$ feuilles. Soit $T$
  un arbre {\`a} $n+1$ feuilles. Pour d{\'e}finir $\psi(T)$, on choisit un
  sommet maximal $s$, comme dans la d{\'e}monstration de la proposition
  \ref{K-binaire}, et on obtient $T$ comme compos{\'e} d'un arbre $t$ {\`a}
  trois feuilles dans un arbre $T'$ {\`a} $n$ feuilles. On pose alors
  $\psi(T)=\psi(T')\circ \psi(t)$. Comme $\phi$ est un morphisme
  d'op{\'e}rades, on a $\phi \psi=\id$.
  
  On montre que $\psi$ ne d{\'e}pend pas des choix faits. On appelle
  secteur d'un arbre plan une paire de feuilles cons{\'e}cutives. Chaque
  choix possible pour obtenir $T$ par composition correspond {\`a} une
  paire de secteurs cons{\'e}cutifs. Consid{\'e}rons deux choix distincts. Si
  les paires de secteurs sont disjointes, l'axiome d'associativit{\'e} des
  op{\'e}rades montre l'ind{\'e}pendance voulue. Sinon, on utilise les
  relations (\ref{relK1}), (\ref{relK6}) ou (\ref{relK7}) de $K'$.
  
  Les autres relations de $K'$ servent {\`a} choisir des repr{\'e}sentants
  particuliers des {\'e}l{\'e}ments de $\quot$ dans $\FF$. Plus pr{\'e}cis{\'e}ment,
  on supprime les compositions de la forme du cot{\'e} gauche de
  (\ref{relK2}), (\ref{relK3}), (\ref{relK4}) et (\ref{relK5}). Les
  repr{\'e}sentants obtenus peuvent s'{\'e}crire comme composition de telle
  fa{\c c}on que leur image par $\phi$ soit une composition du type de
  celles utilis{\'e}es pour d{\'e}finir $\psi$.
  
  Ceci permet de d{\'e}montrer, en utilisant le fait que $\psi$ ne d{\'e}pend
  pas des choix, que $\psi \phi =\id$, donc $\psi$ est un
  inverse de $\phi$ pour les arbres {\`a} $n+1$ feuilles, ce qui termine
  la r{\'e}currence.
\end{proof}

\medskip

La diff{\'e}rentielle de $K'$ est donc uniquement d{\'e}termin{\'e}e par le respect
de la structure d'op{\'e}rade et les conditions suivantes :
\begin{align}
  d(1 < 2)&=-d(1 > 2)=(1 | 2),\\
  d(1 | 2)&=0.
\end{align}

On en d{\'e}duit la description suivante des $K'$-alg{\`e}bres, qui est aussi
celle des $K$-alg{\`e}bres.

 \smallskip

\begin{defi}
  Une $K$-alg{\`e}bre est la donn{\'e}e d'un complexe de cocha{\^\i}nes $V$ et
  d'applications $\succ \,: V\otimes V\to V$, $\prec \,: V\otimes V\to
  V$ et $\times : V\otimes V\to V[-1]$ telles que pour tous $x,y,z$
  dans $V$,
  \begin{align} 
  x\succ(y\prec z)&=(x\succ y)\prec z,\\
  x\succ(y\succ z)&=(x\succ y)\succ z+(x\prec y)\succ z,\\
  (x\prec y)\prec z&=x\prec(y\prec z)+x\prec(y\succ z),\\
  (x\succ y)\times z&=x\succ (y\times z),\\
  x\times (y\prec z)&=(x\times y)\prec z,\\
  (x\prec y)\times z&=(-1)^{y}x\times (y\succ z),\\
   (x\times y)\times z&=x\times (y\times z) 
  \end{align}
  et 
  \begin{align}
    (dx\prec y)+(-1)^{x}(x\prec d y)-d(x\prec y)&=(-1)^{x} x \times y,\\
    (dx\succ y)+(-1)^{x}(x\succ d y)-d(x\succ y)&=(-1)^{x+1}x \times y,\\
    (dx\times y)+(-1)^{x+1}(x\times d y)-d(x\times y)&=0.
   \end{align}
\end{defi}


\smallskip

\begin{rema}
La koszulit{\'e} d'une version non-gradu{\'e}e de l'op{\'e}rade $K'$ est annonc{\'e}e
dans \cite{loday-annonce}. Il est probable que $K'$ est {\'e}galement de
Koszul.
\end{rema}

\smallskip

On observe que le quotient de $K'$ par les {\'e}l{\'e}ments de dimension non
nulle forme une op{\'e}rade (de diff{\'e}rentielle nulle) isomorphe {\`a}
l'op{\'e}rade $\dend$ des alg{\`e}bres dendriformes.

\smallskip

D'autre part, on a un morphisme de $\as$ muni de la diff{\'e}rentielle
nulle dans $K'$.  C'est un quasi-isomorphisme, car les complexes
sous-jacents {\`a} $K'$ sont les complexes cellulaires des polytopes de
Stasheff.

\section{Op{\'e}rades duales}

On donne ici une pr{\'e}sentation par g{\'e}n{\'e}rateurs et relations des
op{\'e}rades duales des trois op{\'e}rades pr{\'e}c{\'e}dentes. On d{\'e}montre aussi que
l'op{\'e}rade (non-sym{\'e}\-trique) duale de $K'$ est isomorphe {\`a} l'op{\'e}rade
$\pasc$ vue comme op{\'e}rade non-sym{\'e}trique.

\smallskip

On obtient les relations du dual d'une op{\'e}rade quadratique binaire par
le calcul de l'orthogonal des relations de cette op{\'e}rade. Plus
pr{\'e}cis{\'e}ment, pour une op{\'e}rade engendr{\'e}e par $E=E(2)$, on consid{\`e}re,
dans l'op{\'e}rade libre $\FF_{E}$ engendr{\'e}e par $E$, le
sous-$\sym_3$-module $R$ de $\FF_{E}(3)$ engendr{\'e} par les
relations. On identifie $\FF_{E}(3)$ avec
\begin{equation*}
E\otimes E\oplus \tau E\otimes E\oplus \tau^2 E\otimes E
\end{equation*}
par l'application 
\begin{equation*}
\tau^k\left( x\otimes y\right) \mapsto \tau^k \left(x \circ_1 y\right).
\end{equation*}
Soit $E'$ le dual de $E$, tensoris{\'e} par la repr{\'e}sentation signe de
$\sym_2$. On d{\'e}finit un couplage entre $\FF_{E}(3)$ et $\FF_{E'}(3)$
en prenant la somme directe de trois couplages entre $E\otimes E$ et
$E'\otimes E'$ pour chacune des puissances de $\tau$. On prend alors
l'orthogonal de $R$ pour ce couplage. De m{\^e}me, la diff{\'e}rentielle de
l'op{\'e}rade duale est obtenue par transposition. On renvoie {\`a}
\cite{ginzburg-kapranov} pour les fondations th{\'e}oriques de cette
dualit{\'e} de Koszul.

\subsection{Trig{\`e}bres}

On travaille dans ce paragraphe avec des op{\'e}rades non-sym{\'e}\-triques.
L'op{\'e}rade $\trias'$ est la duale de Koszul de l'op{\'e}rade $K'$ des
associa{\`e}dres. On note $\trias$ l'op{\'e}rade sym{\'e}trique associ{\'e}e.

L'op{\'e}rade $\trias'$ est engendr{\'e}e par les {\'e}l{\'e}ments $1\dashv2$,
$1\vdash 2$ de dimension $0$ et $1\times 2$ de dimension $1$ soumis
aux relations suivantes :
\begin{align}
  \label{relT1}
  (1 \dashv 2)\circ_1 (1 \vdash 2)&=(1 \vdash 2)\circ_2 (1 \dashv 2),\\
  \label{relT2}
  (1 \dashv 2)\circ_2 (1 \dashv 2)&=(1 \dashv 2)\circ_2 (1 \vdash 2),\\
  \label{relT3}
  (1 \dashv 2)\circ_1 (1 \dashv 2)&=(1 \dashv 2)\circ_2 (1 \vdash 2),\\
  \label{relT4}
  (1 \vdash 2)\circ_1 (1 \vdash 2)&=(1 \vdash 2)\circ_1 (1 \dashv 2),\\
  \label{relT5}
  (1 \vdash 2)\circ_2 (1 \vdash 2)&=(1 \vdash 2)\circ_1 (1 \dashv 2),\\
  \label{relT6}
  (1\times 2)\circ_1 (1 \vdash 2)&=(1 \vdash 2)\circ_2 (1\times 2),\\
  \label{relT7}
  (1\times 2)\circ_2 (1 \dashv 2)&=(1 \dashv 2)\circ_1 (1\times 2),\\
  \label{relT8}
  (1\times 2)\circ_1 (1 \dashv 2)&=(1\times 2)\circ_2 (1 \vdash 2),\\
  \label{relT9}
  (1 \vdash 2)\circ_1 (1\times 2)&=0,\\
  \label{relT10}
  (1\dashv 2)\circ_2 (1 \times 2)&=0,\\
  \label{relT11}
  (1\times 2)\circ_1 (1\times 2)&=-(1\times 2)\circ_2 (1\times 2).
\end{align}

La diff{\'e}rentielle est uniquement d{\'e}termin{\'e}e par le respect de la
structure d'op{\'e}rade et les conditions suivantes :
\begin{align}
  d(1\times 2)&=1\dashv 2- 1\vdash 2,\\
  d(1\vdash 2)&=d(1\dashv 2)=0.
\end{align}

On en d{\'e}duit la description suivante des trig{\`e}bres.

\smallskip

\begin{defi}
  Une trig{\`e}bre est la donn{\'e}e d'un complexe de cha{\^\i}nes $V$,
  d'applications $\vdash\,: V\otimes V\to V$, $\dashv\, : V\otimes
  V\to V $ et $\times : V\otimes V \to V[-1]$ telles que
\begin{align}
  (x\vdash y)\dashv z&=x\vdash (y\dashv z),\\
  x\dashv (y\dashv z)&=x\dashv (y\vdash z)=(x\dashv y)\dashv z,\\
  (x\vdash y)\vdash z&=(x\dashv y)\vdash z=(x\vdash y)\vdash z,\\
  (x\vdash y)\times z&=x\vdash (y\times z),\\
  (x\times y)\dashv z&=x\times (y\dashv z),\\
  (x\dashv y)\times z&=(-1)^{y}x\times (y\vdash z),\\
  (x\times y)\vdash z&=0,\\
  x\dashv (y\times z)&=0,\\
  (x\times y)\times z&=x\times (y\times z)
\end{align}
  et 
  \begin{align}
    (d x\times y)+(-1)^{x+1}(x\times d y) 
    &-d (x\times y)=x \dashv y-    x \vdash y,\\
    (dx \vdash y)+(-1)^{x} (x\vdash dy)
    &-d(x\vdash y)=0,\\
    (dx \dashv y)+(-1)^{x} (x\dashv dy)
    &-d(x\dashv y)=0.
  \end{align}
\end{defi}

On peut donner une description globale de $\trias'$. Si $\PP$ est une
op{\'e}rade sym{\'e}trique, on note $\PP_{\NS}$ l'op{\'e}rade non sym{\'e}trique
obtenue par oubli des actions des groupes sym{\'e}triques.

\begin{prop}
  L'op{\'e}rade non-sym{\'e}trique $\trias'$ est isomorphe {\`a} l'op{\'e}rade
  non-sym{\'e}trique $\pasc_{\NS}$.
\end{prop}
\begin{proof}
  On change de notation : les {\'e}l{\'e}ments $e_1$, $e_2$ et $e_1\wedge e_2$
  de $\pasc$ deviennent les g{\'e}n{\'e}rateurs $1\dashv 2$, $1\vdash 2$ et
  $1\times2$ de $\trias'$. On v{\'e}rifie sans difficult{\'e} que les
  relations qui d{\'e}finissent $\trias$ sont v{\'e}rifi{\'e}es dans
  $\pasc_{\NS}$. Il reste {\`a} d{\'e}montrer que l'on obtient bien ainsi une
  pr{\'e}sentation par g{\'e}n{\'e}rateurs et relations de $\pasc_{\NS}$. On voit
  ais{\'e}ment que les trois {\'e}l{\'e}ments $1\dashv 2$, $1\vdash 2$ et
  $1\times2$ engendrent $\pasc_{\NS}$. Il reste {\`a} d{\'e}montrer que les
  relations de $\trias'$ sont exactement celles de $\pasc_{\NS}$. On ne
  d{\'e}taille pas la preuve, qui est assez semblable {\`a} celle de la
  proposition \ref{pasc-pres}.
\end{proof}

\begin{rema}
En particulier, on obtient ainsi la s{\'e}rie g{\'e}n{\'e}ratrice de $\trias$.
\end{rema}

\smallskip

On observe que la sous-op{\'e}rade (de diff{\'e}rentielle nulle) form{\'e}e par
les {\'e}l{\'e}ments de dimension nulle de $\trias$ est isomorphe {\`a} l'op{\'e}rade
$\dias$ des dig{\`e}bres.

\smallskip

D'autre part, on a un morphisme de $\trias$ dans $\as$ muni de la
diff{\'e}rentielle nulle. C'est un quasi-isomor\-phisme, comme application
duale du quasi-iso\-mor\-phisme de $\as$ dans $K$.

\smallskip

Si $\PP$ est une op{\'e}rade de $\ZZ$-modules, on note $\PP_{\QQ}$
l'op{\'e}rade d'espace vectoriels sur $\QQ$ obtenue par extension des
scalaires.

\begin{prop}
  On a un isomorphisme d'op{\'e}rades sym{\'e}triques
\begin{equation*}	
\left( \pasc\otimes\as\right)_{\QQ}\simeq \trias_{\QQ}.
\end{equation*}
\end{prop}
\begin{proof}
  On v{\'e}rifie que les relations de $\trias$ sont celles de l'op{\'e}rade
  $\pasc \circ \as$ o{\`u} le produit $\circ$ est le produit de Manin des
  op{\'e}rades quadratiques binaires, voir \cite[2.2]{ginzburg-kapranov}.
  Par cons{\'e}quent, $\trias$ est la sous-op{\'e}rade de $\pasc\otimes \as$
  engendr{\'e}e par $\pasc(2)\otimes \as(2)$. Comme on a {\'e}galit{\'e} des
  s{\'e}ries g{\'e}n{\'e}ratrices, on peut en d{\'e}duire un isomorphisme des op{\'e}rades
  tensoris{\'e}es par $\QQ$.
\end{proof}

\subsection{Alg{\`e}bres $\pidual$}

Cette op{\'e}rade est la duale de Koszul de l'op{\'e}rade $\Pi$ sur les
permuto{\`e}dres. Elle est engendr{\'e}e par les {\'e}l{\'e}ments suivants : $[1,2]$
sym{\'e}trique de dimension $1$, $\langle 1, 2\rangle$ et $\langle 2,
1\rangle$ de dimension $0$, modulo l'id{\'e}al engendr{\'e} par les relations
suivantes :
\begin{align}
  \label{relU1}
  \langle 1,2 \rangle \circ_1 \langle 1,2 \rangle - \tau \langle 1,2
  \rangle \circ_1 \langle 2,1 \rangle
  &+\tau^2 \langle 2,1 \rangle \circ_1 \langle 2,1 \rangle,\\
  \label{relU2}
  \tau^2 [1,2]\circ_1 \langle 1,2 \rangle +\tau [1,2]\circ_1 \langle
  2,1 \rangle
  &=\langle 1,2 \rangle \circ_1 [1,2],\\
  \label{relU3}
  \langle 2,1 \rangle \circ_1 [1,2]&= 0,\\
  \label{relU4}
  [1,2]\circ_1 [1,2]+\tau[1,2]\circ_1 [1,2]&+\tau^2[1,2]\circ_1 [1,2]=0,
\end{align}
o{\`u} $\tau$ est le cycle $3\to 2\to 1\to 3$. La diff{\'e}rentielle est
uniquement d{\'e}termin{\'e}e par le respect de la structure d'op{\'e}rade et les
conditions suivantes :
\begin{align}
  d([1,2])&= \langle 1,2 \rangle+\langle 2,1 \rangle,\\
  d(\langle 1,2 \rangle)&=0.
\end{align}

On en d{\'e}duit la description suivante des $\pidual$-alg{\`e}bres. 
\smallskip

\begin{defi}
  Une $\pidual$-alg{\`e}bre est la donn{\'e}e d'un complexe de cha{\^\i}nes $V$,
  d'appli\-cations $\langle\cdot\, , \,\cdot \rangle : V\otimes V\to V $
  et $[\cdot\, , \, \cdot] : V\otimes V\to V[-1]$ telles que
\begin{align}
  \langle x,\langle y,z \rangle \rangle &= \langle \langle x,y
  \rangle,z \rangle   - \langle \langle x,z \rangle ,y\rangle  ,\\
  \langle[x,y], z\rangle =(-1)^{yz+z}&[\langle x,z\rangle,y]
  +[x,\langle y,z\rangle],\\
  \langle x,[y, z]\rangle &=0,\\
  [x,y]&=(-1)^{xy+x+y}[y,x],\\
  (-1)^{y+zx}[[x,y],z]+(-1)^{z+xy}&[[y,z],x]+(-1)^{x+yz}[[z,x],y]=0
\end{align}
   et 
  \begin{align}
    [dx, y]+(-1)^{x+1}[x, dy]- d([x, y])
    &=\langle x, y\rangle + (-1)^{xy}\langle y ,x\rangle,\\
    (\langle dx,y\rangle)+(-1)^{x}(\langle x,dy\rangle &-d(\langle x,y
    \rangle)=0.
  \end{align}
\end{defi}

La sous-op{\'e}rade (de diff{\'e}rentielle nulle) form{\'e}e par les {\'e}l{\'e}ments de
dimension nulle de $\pidual$ est isomorphe {\`a} l'op{\'e}rade $\leib$ des
alg{\`e}bres de Leibniz.

\smallskip

D'autre part, on a un morphisme de $\pidual$ dans $\lie$ muni de la
diff{\'e}rentielle nulle. C'est un quasi-isomorphisme, comme application
duale du quasi-isomorphisme de $\com$ dans $\Pi$.

\begin{prop}
  On a un isomorphisme $\left( \pasc\otimes\lie\right) _{\QQ}\simeq
  \pidual_{\QQ}$.
\end{prop}
\begin{proof}
  On v{\'e}rifie que les relations de $\pidual$ sont celles de l'op{\'e}rade
  $\pasc \circ \lie$ o{\`u} le produit $\circ$ est le produit de Manin des
  op{\'e}rades quadratiques binaires, voir \cite[2.2]{ginzburg-kapranov}.
  Par cons{\'e}quent, $\pidual$ est la sous-op{\'e}rade de $\pasc\otimes \lie$
  engendr{\'e}e par $\pasc(2)\otimes \lie(2)$. Comme on a {\'e}galit{\'e} des
  s{\'e}ries g{\'e}n{\'e}ratrices, on peut en d{\'e}duire un isomorphisme des op{\'e}rades
  tensoris{\'e}es par $\QQ$.
\end{proof}

\subsection{Alg{\`e}bres $\Lambda$}

Cette op{\'e}rade est la duale de Koszul de l'op{\'e}rade $\pasc$ sur les
simplexes. Elle a une saveur ``espace de configurations'', voir par
exemple sa s{\'e}rie caract{\'e}ristique suppos{\'e}e dans le tableau
r{\'e}capitulatif. Elle est engendr{\'e}e par les {\'e}l{\'e}ments suivants : $[1,2]$
sym{\'e}trique de dimension $1$, $1\courb 2$ et $2\courb 1$ de dimension
$0$, modulo l'id{\'e}al engendr{\'e} par les relations suivantes :
\begin{align}
  \label{relM1}
  \tau (2 \courb 1)\circ_1 (1 \courb 2)-(1 \courb 2)\circ_1& (1 \courb 2)\\
\notag  =\tau (2 \courb 1)\circ_1& (2 \courb 1)
  -\tau^2(1 \courb 2)\circ_1 (2 \courb 1),\\
  \label{relM2}
  [1,2]\circ(1\courb 2)-\tau[1,2]\circ(2\courb 1)
  &=\tau^2 (1 \courb 2)\circ_1 [1,2],\\
  \label{relM3}
  [1,2]\circ_1 [1,2]+\tau[1,2]\circ_1 [1,2]&+\tau^2[1,2]\circ_1 [1,2]=0,
\end{align}
o{\`u} $\tau$ est le cycle $3\to 2\to 1\to 3$. La diff{\'e}rentielle est
uniquement d{\'e}termin{\'e}e par le respect de la structure d'op{\'e}rade et les
conditions suivantes :
\begin{align}
  d(1\courb 2)&=[1,2],\\
  d([1,2])&=0.
\end{align}

On en d{\'e}duit la description suivante des $\Lambda$-alg{\`e}bres. 
\smallskip

\begin{defi}
  Une $\Lambda$-alg{\`e}bre est la donn{\'e}e d'un complexe de cocha{\^\i}nes $V$,
  d'appli\-cations $\circ : V\otimes V\to V $ et $[\cdot\, , \, \cdot] :
  V\otimes V\to V[-1] $ telles que
  \begin{align}
    (x\courb y)\courb z-x\courb(y\courb z)&=
    (x\courb z)\courb y-x\courb(z\courb y),\\
    [x,y]\courb z&=(-1)^{zy+z}[x\courb z,y]+[x,y\courb z],\\
    [x,y]&=(-1)^{xy+x+y}[y,x],\\
    (-1)^{y+zx}[[x,y],z]+(-1)^{z+xy}&[[y,z],x]+(-1)^{x+yz}[[z,x],y]=0
  \end{align}
  et
  \begin{align}
    (dx\courb y)+(-1)^{x}(x\courb d y)&-d(x\courb y)=(-1)^{x}[x,y],\\
    ([dx,y])+(-1)^{x+1}[x,dy]&-d([x,y])=0.
  \end{align}
\end{defi}

On observe que le quotient de $\Lambda$ par les {\'e}l{\'e}ments de dimension
non nulle forme une op{\'e}rade (de diff{\'e}rentielle nulle) isomorphe {\`a}
l'op{\'e}rade $\prelie$ des alg{\`e}bres pr{\'e}-Lie.

\smallskip

D'autre part, on a un morphisme de l'op{\'e}rade $\lie$ (munie de la
diff{\'e}rentielle nulle) dans $\Lambda$. C'est un quasi-isomorphisme
comme application duale du quasi-isomorphisme de $\pasc$ dans $\com$.

\section{Morphismes horizontaux}

Il reste, pour compl{\'e}ter la description du diagramme (\ref{grosdiag}),
{\`a} d{\'e}crire les applications horizontales des lignes $2$ et $4$. Les
autres morphismes horizontaux ont d{\'e}j{\`a} {\'e}t{\'e} d{\'e}finis, voir par exemple
\cite{dialgebra} ; on peut les retrouver par restriction ou passage au
quotient de ceux qu'on d{\'e}finit ci-dessous. On consid{\`e}re ici $K$ et
$\trias$ comme des op{\'e}rades avec action libre des groupes sym{\'e}triques.

\smallskip

On a un morphisme d'op{\'e}rades de $K$ dans $\Pi$ d{\'e}fini au niveau des
alg{\`e}bres par les formules suivantes :
\begin{align}
  x\prec y&:=x\prec y,\\
  x\succ y&:=(-1)^{xy}y\prec x,\\
  x\times y&:=x\times y.
\end{align}

\smallskip

Le morphisme d'op{\'e}rades de $\Lambda$ dans $K$ est d{\'e}fini au niveau des
alg{\`e}bres par les formules suivantes :
\begin{align}
  x\courb y&:=x\prec y-(-1)^{xy} y\succ x,\\ 
  [x,y]&:=x\times y+(-1)^{xy+x+y} y\times x. 
\end{align}

\smallskip

Les formules suivantes d{\'e}finissent, au niveau des alg{\`e}bres, un
morphisme d'op{\'e}rades de $\trias$ dans $\pasc$ :
\begin{align}
    x \dashv y&:= x \dashv y,\\
    x \vdash y&:=(-1)^{xy}y\dashv x,\\
    x \times y&:=x\times y.
\end{align}

\smallskip
On a un morphisme d'op{\'e}rades de $\pidual$ dans $\trias$ d{\'e}fini au niveau des
alg{\`e}bres par les formules suivantes :
\begin{align}
  \langle x,y \rangle&:=x\dashv y -(-1)^{xy} y \vdash x,\\
  [x,y]&:=x\times y+(-1)^{xy+x+y} y\times x. 
\end{align}

\smallskip

On v{\'e}rifie alors sans difficult{\'e} que

\begin{prop}
  Le diagramme (\ref{grosdiag}) est form{\'e} de carr{\'e}s commutatifs.
\end{prop}

De plus, chacune des lignes de ce diagramme est ``exacte'', au sens
suivant : l'op{\'e}rade de gauche est le quotient de l'op{\'e}rade du milieu
par l'id{\'e}al engendr{\'e} par l'image de l'id{\'e}al d'augmentation de
l'op{\'e}rade de droite.

\section{S{\'e}ries g{\'e}n{\'e}ratrices}

On utilise la convention suivante pour les s{\'e}ries g{\'e}n{\'e}ratrices
d'op{\'e}rades diff{\'e}ren\-tielles gradu{\'e}es :

\begin{equation*}
  g_{\PP}(x,t)=\sum_{n\geq 1}\sum_{k}\dim \PP^k(n) (-t)^k \frac{x^n}{n!}.
\end{equation*}

La s{\'e}rie g{\'e}n{\'e}ratrice de la suspension $\Sigma \PP$ d'une op{\'e}rade $\PP$
v{\'e}rifie
\begin{equation*}
  g_{\Sigma\PP}(x,t)=-g_{\PP}(-tx,t)/t.
\end{equation*}

Les s{\'e}ries g{\'e}n{\'e}ratrices des op{\'e}rades qui apparaissent dans le
diagramme (\ref{grosdiag}) sont r{\'e}capitul{\'e}es dans le tableau
ci-dessous. Celles des op{\'e}rades des premi{\`e}re, troisi{\`e}me et cinqui{\`e}me
lignes sont bien connues. Celles des op{\'e}rades $\pasc$, $K$, $\Pi$ et
$\trias$ s'obtiennent ais{\'e}ment en partant de la description explicite
de ces op{\'e}rades. Celle de $\pidual$ r{\'e}sulte de la koszulit{\'e} de $\Pi$.
Enfin, celle de $\Lambda$ s'obtient par inversion de celle de $\pasc$,
sous r{\'e}serve que les op{\'e}rades $\pasc$ et $\Lambda$ soient de Koszul.

\bigskip

\begin{tabular}{|p{2.5cm}|p{4cm}|p{4cm}|}
  \hline 
  $\zin$ : & $\dend$ : & $\prelie$ : \\  
  $\frac{x}{1-x}$ & $\frac{(1-2 x)-\sqrt{(1-4 x)}}{2 x}$ & \\
  $\sum n! \frac{x^n}{n!}$ & $\sum \frac{1}{n+1}\binom{2 n}{n} x^n$ 
  & $\sum n^{n-1} \frac{x^n}{n!}$ \\ 
  \hline 
  $\Pi$ : & $K$ : & $\Lambda$ : \\ 
  $\frac{e^{tx}-1}{1+(t-1)e^{tx}}$ & 
  $\frac{(1-2 x+x t)-\sqrt{1-4 x+ 2 x t +x^2 t^2}}{2 x (1-t)}$ &  \\
  $\sum \pi_{n,k} (-t)^k \frac{x^n}{n!}$ & $\sum c_{n,k} (-t)^k x^n $ 
  & $\sum \prod_{k=1}^{n-1}(n-kt)\frac{x^n}{n!}$\\
  \hline $\com$ : & $\as$ : & $\lie$ : \\ 
  $e^{x}-1$ & $\frac{x}{1-x} $ & $-\ln(1-x)$ \\
  $\sum \frac{x^n}{n!}$ & $\sum x^n$ & $\sum (n-1)! \frac{x^n}{n!}$ \\
  \hline $\pasc$ : & $\trias$ : & $\pidual$ : \\ 
  $\frac{e^x-e^{(1-t)x}}{t}$ &
  $\frac{1}{t}(\frac{x}{1-x}-\frac{(1-t)x}{1+(t-1)x})$ & 
  $\frac{-\ln(1-x)+\ln(1+(t-1)x)}{t}$ \\
  $\sum \binom{n}{k+1} (-t)^k \frac{x^n}{n!}$ &
  $\sum \binom{n}{k+1} (-t)^k x^n$ & 
  $\sum \binom{n}{k+1} (n-1)! (-t)^k \frac{x^n}{n!}$ \\
  \hline $\perm$ : & $\dias$ : & $\leib$ : \\ 
  $x e^{x}$ & $\frac{x}{(1-x)^2}$ & $\frac{x}{1-x} $ \\
  $\sum n \frac{x^n}{n!}$ & $\sum n\, x^n$ & $\sum n! \frac{x^n}{n!}$ \\
  \hline
\end{tabular}

\smallskip

Dans ce tableau, $\pi_{n,k}$ (resp. $c_{n,k}$) d{\'e}signe le nombre de
faces de dimension $k$ du permuto{\`e}dre (resp. de l'associa{\`e}dre) de
dimension $n-1$.

\providecommand{\bysame}{\leavevmode ---\ }
\providecommand{\og}{``}
\providecommand{\fg}{''}
\providecommand{\smfandname}{et}
\providecommand{\smfedsname}{{\'e}ds.}
\providecommand{\smfedname}{{\'e}d.}
\providecommand{\smfmastersthesisname}{M{\'e}moire}
\providecommand{\smfphdthesisname}{Th{\`e}se}

\end{document}